\documentclass[12pt,leqno]{article}
\usepackage{amsfonts}
\usepackage{amsmath, amsthm, amsfonts, amssymb, color}
\usepackage{mathrsfs}
\theoremstyle{definition}

\newcommand{\scr}[1]{\mathscr #1}
\definecolor{wco}{rgb}{0.5,0.2,0.3}
\allowdisplaybreaks[2]
\numberwithin{equation}{section} \theoremstyle{remark}

\newcommand{\ua}{\uparrow}

\title{
{\bf  Weighted Super Poincar\'e Inequalities for Infinite-Dimensional Extension of the Dirichlet Distribution}
\footnote{corresponding author email address: 201321130173@mail.bnu.edu.cn}}
\author{
{\bf Weiwei Zhang}\\
\footnotesize{School of Mathematical Sciences, Beijing Normal University, Beijing 100875, China}
}

\begin{document}
\def\tttext#1{{\normalfont\ttfamily#1}}
\def\R{\mathbb R}  \def\ff{\frac} \def\ss{\sqrt} \def\B{\mathbf B}
\def\N{\mathbb N} \def\kk{\kappa} \def\m{{\bf m}}
\def\dd{\delta} \def\DD{\Delta} \def\vv{\varepsilon} \def\rr{\rho}
\def\<{\langle} \def\>{\rangle} \def\GG{\Gamma} \def\gg{\gamma}
  \def\nn{\nabla} \def\pp{\partial} \def\EE{\scr E}
\def\d{\text{\rm{d}}} \def\bb{\beta} \def\aa{\alpha} \def\D{\scr D}
  \def\si{\sigma} \def\ess{\text{\rm{ess}}}
\def\beg{\begin} \def\beq{\begin{equation}}  \def\F{\scr F}
\def\Ric{\text{\rm{Ric}}} \def\Hess{\text{\rm{Hess}}}
\def\e{\text{\rm{e}}} \def\ua{\underline a} \def\OO{\Omega}  \def\oo{\omega}
 \def\tt{\tilde} \def\Ric{\text{\rm{Ric}}}
\def\cut{\text{\rm{cut}}} \def\P{\mathbb P}
\def\C{\scr C}     \def\E{\mathbb E}\def\y{{\bf y}}
\def\Z{\mathbb Z} \def\II{\mathbb I}
  \def\Q{\mathbb Q}  \def\LL{\Lambda}\def\L{\scr L}
  \def\B{\scr B}    \def\ll{\lambda} \def\a{{\bf a}} \def\b{{\bf b}}
\def\vp{\varphi}\def\H{\mathbb H}\def\ee{\mathbf e}\def\x{{\bf x}}
\def\gap{{\rm gap}}\def\PP{\scr P}\def\p{{\mathbf p}}\def\NN{\mathbb N}
\def\cA{\scr A} \def\cQ{\scr Q}\def\cK{\scr K}
\def\LS{C_{LS}}
\maketitle
\begin{abstract}  For the infinite-dimensional extension of the Dirichlet distribution, the super Poincar\'e inequality does not hold based on the result in \cite{WZ}, so we establish the weighted super Poincar\'e inequalities for this measure with respect to two different Dirichlet forms respectively.

 \end{abstract} \noindent
   \vskip 2cm

\section{Introduction}
Let $ n\ge 1$ be a natural number, and let  $ \aa=(\aa_1,\cdots, \aa_{n+1})\in (0,\infty)^{n+1}.$   The Dirichlet distribution $\mu^{(n)}_{\aa}$ with parameter $\aa$ is a probability measure on the set
 $$\DD^{(n)}:= \bigg\{ x=(x_i)_{1\le i\le n}\in [0,1]^n:\ |x|_1:=\sum_{i=1}^n x_i\le 1\bigg\}$$ with density function
\beq\label{RR} \rr(x):= \ff{\GG(|\aa|_1)}{\prod_{1\le i\le n+1} \GG(\aa_i)} (1-|x|_1)^{\aa_{n+1}-1}\prod_{1\le i\le n} x_i^{\aa_i-1},\ \ x=(x_i)_{1\le i\le n}\in \DD^{(n)},\end{equation}
where $|\aa|_1:=\sum_{i=1}^{n+1}\aa_i$, denoted by $D(\alpha_{1},\alpha_{2},\cdots,\alpha_{n},\alpha_{n+1}).$
This   distribution  arises naturally in Bayesian inference as conjugate prior  for categorical distribution, and it describes   the distribution of  allelic frequencies in  population genetics,  see for instance \cite{P1}.
We introduce two generators of diffusion processes whose stationary distributions both are the Dirichlet distribution.
\beg{enumerate} \item[(1)] The first type operator and Dirichlet form:
\beg{equation}\beg{split}\label{one}&L_{\alpha,1}^{(n)}f(x)
=\sum_{i=1}^{n}x_{i}(1-|x|_{1})(\partial_{i}^{2}f)(x)+\{\alpha_{i}(1-|x|_{1})
-\alpha_{\infty}x_{i}\}\partial_{i}f(x),\quad f\in C^{2}(\mathbb{R}^{n}).\\
&\EE_{\alpha,1}^{(n)}(f,g):=\mu_{\alpha,\alpha_{\infty}}^{(n)}\bigg((1-|x|_{1})\sum_{i=1}^{n}x_{i}\partial_{i}f\partial_{n}g\bigg),\quad f, g\in C^{1}(\mathbb{R}^{n}).\end{split}\end{equation}
\item[(2)]The second type operator and Dirichlet form:

\beg{equation}\beg{split}&L_{\alpha,2}^{(n)}f(x)
=\sum_{i,j=1}^{n}x_{i}(\delta_{ij}-x_{j})(\partial_{ij}f)(x)+\sum_{i=1}^{n}(\alpha_{i}-|\alpha|_{1}x_{i})\partial_{i}f(x)
,\quad f\in C^{2}(\mathbb{R}^{n}).\\
&\EE_{\alpha,2}^{(n)}(f,g):=\mu_{\alpha,\alpha_{\infty}}^{(n)}\bigg(\sum_{i,j=1}^{n}x_{i}(\delta_{ij}-x_{j})\partial_{i}f\partial_{j}g\bigg),\quad f,g\in C^{1}(\mathbb{R}^{n}).\end{split}\end{equation}
\end{enumerate}
\cite{FMW} extend the measure $\mu^{(n)}_{\aa}$ and the operator \eqref{one} to the infinite-dimensional case.
Consider the infinite-dimensional simplex
$$\DD^{(\infty)}:= \bigg\{ x\in [0,1]^{\mathbb{N}}:\ |x|_1:=\sum_{i=1}^\infty x_i\le 1\bigg\},  $$
which is equipped with the $L^{1}$-metric $|x-y|_{1}.$
Let
$$\alpha^{(n)}:=(\alpha_{1},\cdots,\alpha_{n-1},\sum_{i\geq n}\alpha_{i},\alpha_{\infty}).$$
We define a probability measure on $ \Delta^{(n)},$
$$\mu^{(n)}_{\alpha^{(n)}}:= D(\alpha_{1},\alpha_{2},\cdots,\alpha_{n-1},\Sigma_{i\geq n}\alpha_{i},\alpha_{\infty});$$
a probability measure on $ \Delta^{(\infty)},$
$$\mu^{(m)}_{\alpha,\alpha_{\infty}}(dx):=\mu^{(n)}_{\alpha^{(n)}}(dx_{1},\cdots,dx_{n})\prod_{i>n}\delta_{0}(dx_{i}).$$
According to \cite{FMW}, we know $\{\mu^{(m)}_{\alpha,\alpha_{\infty}}\}_{m\geq1}$ weakly converges and we denote the limit as  $\mu_{\alpha,\alpha_{\infty}}^{(\infty)}.$
$\mathcal{F}C^{p}$ is the $C^{p}$-cylindrical functions for $p\geq 1,$
$$\mathcal{F}C^{p}:=\{\Delta^{(\infty)}\ni x:=(x_{i})_{i\geq 1}\mapsto f(x_{1},\cdots,x_{m}):m\geq 1, f\in C^{p}(\Delta^{(m)})\}.$$

 For the Dirichlet distribution and the first type Dirichlet form, S. Feng, L. Miclo and F.-Y. Wang have established the Poincar\'e inequalities for the finite-dimensional case and the infinite-dimensional case in \cite{FMW}. Then the super Poincar\'e inequality for the finite-dimensional case is established in \cite{WZ} by F.-Y. Wang and the author of this paper. What's more, \cite{FMW} also proved that the form
\beg{align*}\EE_{\alpha,1}^{(\infty)}(f,f):=\mu_{\alpha,\alpha_{\infty}}^{(\infty)}\bigg((1-|x|_{1})\sum_{n=1}^{\infty}x_{n}\partial_{n}f\partial_{n}g\bigg)\quad f,g\in\mathcal{F}C^{1}\end{align*}
is closable in $L^{2}(\mu_{\alpha,\alpha_{\infty}}^{(\infty)})$ and the closure is a symmetric Dirichlet form. The generator\\ $(L_{\alpha,1}^{(\infty)},\D(L_{\alpha,1}^{(\infty)}))$ of the Dirichlet form $\EE_{\alpha,1}^{(\infty)}$ satisfies $ \mathcal{F}C^{2}\subset\D(L_{\alpha,1}^{(\infty)}),$
\beg{align*}L_{\alpha,1}^{(\infty)}f(x)
=\sum_{n=1}^{\infty}\bigg(x_{n}(1-|x|_{1})\partial_{n}^{2}f(x)+\{\alpha_{n}(1-|x|_{1})
-\alpha_{\infty}x_{n}\}\partial_{n}f(x)\bigg),\quad f\in\mathcal{F}C^{2}.\end{align*}

In \cite{S}, Stannat have established the Poincar\'e inequality for the Dirichlet distribution with respect to the second type Dirichlet form and also established the Poincar\'e inequality for the Fleming-Viot process. In \cite{WZ}, F.-Y. Wang and the author of this paper have established the super Poincar\'e inequality for the second type  Dirichlet form about finite-dimensional case. Following the idea of \cite{FMW}, we can also prove that the form
\beg{align*}\EE_{\alpha,2}^{(\infty)}:=\mu_{\alpha,\alpha_{\infty}}^{(\infty)}\bigg(\sum_{i,j=1}^{\infty}x_{i}(\delta_{ij}-x_{j})\partial_{i}f\partial_{j}g\bigg)\quad f,g\in\mathcal{F}C^{1}\end{align*}
is closable in $L^{2}(\mu_{\alpha,\alpha_{\infty}}^{(\infty)})$ and the closure is a symmetric Dirichlet form. The generator\\ $(L_{\alpha,2}^{(\infty)},\D(L_{\alpha,2}^{(\infty)}))$ of the Dirichlet form $\EE_{\alpha,2}^{(\infty)}$ satisfies $ \mathcal{F}C^{2}\subset\D(L_{\alpha,2}^{(\infty)}),$
\beg{align*}L_{\alpha,2}^{(\infty)}f(x)
=\sum_{i,j=1}^{\infty}x_{i}(\delta_{ij}-x_{j})(\partial_{ij}f)(x)+\sum_{i=1}^{\infty}(\alpha_{i}-|\alpha|_{1}x_{i})\partial_{i}f(x),\quad f\in\mathcal{F}C^{2}.\end{align*}

In \cite{S}, Stannat proved the invalidity of the log-Sobolev inequality when the state space contains infinite states. In fact, by the same method, we can prove the F-Sobolev inequality doesn't hold, then neither does the super Poincar\'e inequality hold.

In this paper, we firstly introduce some properties of the Dirichlet distribution in section one, then we establish the weighted super Poincar\'e inequality for the Dirichlet form $\EE_{\alpha,1}^{(\infty)} $ and the measure  $\mu_{\alpha,\alpha_{\infty}}^{(\infty)}$ in section two.
At last, we obtain the weighted super Poincar\'e inequality for the Dirichlet form $\EE_{\alpha,2}^{(\infty)} $ and the measure $\mu_{\alpha,\alpha_{\infty}}^{(\infty)}$ in section three.

\section{Property of the Dirichlet distribution}
The Dirichlet distribution possesses many nice properties. One of the properties is the partition property of $\mu_{\alpha}^{(n)}$ for $\alpha\in(0,\infty)^{n}.$
Let $(X_{1}, \cdots, X_{n})$ have law $\mu_{\alpha}^{(n)},$ let $\{A_{1}, \cdots, A_{k}\}$ be a partition of the set $\{1, \cdots, n\},$ and set
$$Y_{j}=\sum_{r\in A_{j}}X_{r},\quad \beta_{j}=\sum_{r\in A_{j}}\alpha_{r},\quad j=1, \cdots, k.$$
Then $(Y_{1}, \cdots, Y_{k})$ has law $\mu_{\beta}^{(k)}$ with parameter $\beta:=(\beta_{1}, \cdots, \beta_{k})\in (0,\infty)^{k}.$

We give another property of the Dirichlet distribution $\mu_{\alpha}^{(n)}$ which we will use later.
For any $m> n,$ we define the map
$$ T_{m}:\Delta^{(n)}\times\Delta^{(m-n)}\rightarrow \Delta^{(m)},$$
$$T_{m}:(x_{1},x_{2},\cdots,x_{m})\longmapsto (x_{1}(1-\Sigma_{n< i\leq m}x_{i}),\cdots,x_{n}(1-\Sigma_{n< i\leq m}x_{i}),x_{n+1},\cdots,x_{m}).$$
$$ T:\Delta^{(n)}\times\Delta^{(\infty)}\rightarrow \Delta^{(\infty)},$$
$$T:(x_{1},x_{2},\cdots)\longmapsto (x_{1}(1-\Sigma_{ i> n}x_{i}),\cdots,x_{n}(1-\Sigma_{i> n}x_{i}),x_{n+1},\cdots).$$
Fixed $n\geq 1,$  $\forall f\in \mathcal{F}C^{1}(\Delta^{(\infty)}),$ we denote the number of the variables of $f$ as $m.$
Let
$$\mu_{1}:=D(\alpha_{1},\alpha_{2},\cdots,\alpha_{n},\alpha_{\infty})$$ is a probability measure on $ \Delta^{(n)}.$
$$\mu_{2}^{(m)}:=D(\alpha_{n+1},\alpha_{n+2},\cdots,\alpha_{m-1},\Sigma_{i\geq m}\alpha_{i},\Sigma_{i=1}^{n}\alpha_{i}+\alpha_{\infty})$$ is a probability measure on $ \Delta^{(m-n)}.$
$$ \mu_{2,\infty}^{(m)}(dx):= \mu_{2}^{(m)}(dx_{1},\cdots,dx_{m})\prod_{i>m}\delta_{0}(dx_{i})$$ is a probability measure on $\Delta^{(\infty)}.$
By the same proof of the weak convergence of $\{\mu^{(m)}_{\alpha,\alpha_{\infty}}\}_{m\geq1}$ in \cite{FMW}, we know $\{\mu_{2}^{(m)}\}_{m\geq1}$ weak converges and we define the limit as $\mu_{2}^{(\infty)}.$
\beg{prp}\label{projection}
$$\mu_{\alpha,\alpha_{\infty}}^{(\infty)}(f)=\mu_{1}(\mu_{2}^{(\infty)}(f\circ T)):=\int_{\Delta^{(n)}}\int_{\Delta^{(\infty)}}f\circ T(x,y)\mu_{2}^{(\infty)}(dy)\mu_{1}(dx).$$
\end{prp}
\beg{proof}
As $f$ has $m$ variables,
so
\beg{align*}&\mu_{1}(\mu_{2,\infty}^{(m)}(f\circ T_{m}))=\mu_{1}(\mu_{2}^{(m)}(f\circ T_{m})) \\
&=\frac{\Gamma(\sum_{i=1}^{n}\alpha_{i}+\alpha_{\infty})}{\prod_{i=1}^{n}\Gamma(\alpha_{i})\Gamma(\alpha_{\infty})}
\frac{\Gamma(\sum_{i=1}^{\infty}\alpha_{i})}{\prod_{i=n+1}^{m-1}\Gamma(\alpha_{i})\Gamma(\sum_{i\geq m}\alpha_{i})\Gamma(\sum_{i=1}^{n}\alpha_{i}+\alpha_{\infty})}\\
&\qquad\cdot\int_{\Delta^{(m)}}f(T_{m}(x))x_1^{\aa_1-1}\cdots x_n^{\aa_n-1}
\cdot\bigg(1-\sum_{i=1}^{n}x_i\bigg)^{\aa_\infty-1}x_{n+1}^{\aa_{n+1}-1}\cdots x_{m-1}^{\aa_{m}-1}x_{m}^{\Sigma_{i\geq m}\alpha_{i}-1}\\
&\qquad\cdot\bigg(1-\sum_{i=n+1}^{m}x_i\bigg)^{\sum_{i=1}^{n}\alpha_{i}+\aa_\infty-1}dx_{1}\cdots dx_{m}\\
&=\frac{\Gamma(\sum_{i=1}^{\infty}\alpha_{i})}{\prod_{i=1}^{m-1}\Gamma(\alpha_{i})\Gamma(\sum_{i\geq m}\alpha_{i})\Gamma(\alpha_{\infty})}\cdot\int_{\Delta^{(m)}}f(T_{m}(x))x_1^{\aa_1-1}\cdots x_n^{\aa_n-1}\bigg(1-\sum_{i=n+1}^{m}x_i\bigg)^{\sum_{i=1}^{n}\alpha_{i}-n} \\
&\qquad\cdot x_{n+1}^{\aa_{n+1}-1}\cdots x_{m-1}^{\aa_{m}-1}\cdot x_{m}^{\sum_{i\geq m}\alpha_{i}-1}
\cdot\bigg[\bigg(1-\sum_{i=n+1}^{m}x_i\bigg)\bigg(1-\sum_{i=1}^{n}x_i\bigg)\bigg]^{\aa_\infty-1}\\
&\qquad\cdot\bigg(1-\sum_{i=n+1}^{m}x_i\bigg)^{n}dx_{1}\cdots dx_{m}\\
&=\frac{\Gamma(\sum_{i=1}^{\infty}\alpha_{i})}{\prod_{i=1}^{m-1}\Gamma(\alpha_{i})\Gamma(\sum_{i\geq m}\alpha_{i})\Gamma(\alpha_{\infty})}\\
&\qquad \cdot\int_{\Delta^{(m)}}f(y)y_1^{\aa_1-1}\cdots y_{m}^{\sum_{i\geq m}\alpha_{i}-1}
\bigg(1-\sum_{i=1}^{m}y_i\bigg)^{\aa_\infty-1}\bigg(1-\sum_{i=n+1}^{m}x_i\bigg)^{n}\\
&\qquad\cdot\det(\nabla T^{-1}_{m}(y)) dy_{1}\cdots dy_{m}\\
&=\frac{\Gamma(\sum_{i=1}^{\infty}\alpha_{i})}{\prod_{i=1}^{m-1}\Gamma(\alpha_{i})\Gamma(\sum_{i\geq m}\alpha_{i})\Gamma(\alpha_{\infty})}\cdot\int_{\Delta^{(m)}}f(y)y_1^{\aa_1-1}\cdots y_{m}^{\sum_{i\geq m}\alpha_{i}-1}\\
&\qquad\cdot\bigg(1-\sum_{i=1}^{m}y_i\bigg)^{\aa_\infty-1}dy_{1}\cdots dy_{m}\\
&=\mu^{(m)}_{\alpha^{(m)}}(f)=\mu^{(m)}_{\alpha,\alpha_{\infty}}(f).\end{align*}
We have used
$$\det(\nabla T^{-1}_{m}(y))=\bigg(1-\sum_{i=n+1}^{m}x_i\bigg)^{-n},$$
\beg{align*}
\bigg[\bigg(1-\sum_{i=n+1}^{m}x_i\bigg)\bigg(1-\sum_{i=1}^{n}x_i\bigg)\bigg]=\bigg(1-\sum_{i=n+1}^{m}x_i\bigg)-\bigg(1-\sum_{i=n+1}^{m}x_i\bigg)\sum_{i=1}^{n}x_i
=1-\sum_{i=1}^{m}y_i.\end{align*}
So we get
\beq\label{step}\mu_{1}(\mu_{2,\infty}^{(m)}(f\circ T_{m}))=\mu^{(m)}_{\alpha,\alpha_{\infty}}(f).\end{equation}
Let $m\rightarrow\infty$ at the both side of \eqref{step},
we get
$$\mu_{\alpha,\alpha_{\infty}}^{(\infty)}(f)=\mu_{1}(\mu_{2}^{(\infty)}(f\circ T)):=\int_{\Delta^{(n)}}\int_{\Delta^{(\infty)}}f\circ T(x,y)\mu_{2}^{(\infty)}(dy)\mu_{1}(dx) .$$
\end{proof}

\section{The weighted super Poincar\'e inequality for the first type Dirichlet form }
According to Theorem 1.1 in \cite{WZ}, from the Nash inequality, we can get the super Poincar\'e inequality for the Dirichlet distribution $\mu^{(n)}_{\aa}$ and Dirichlet form $\EE_{\alpha,1}^{(n)}.$ Combining the super Poincar\'e inequality with the Poincar\'e inequality, we can obtain the weighted super Poincar\'e inequality for infinite-dimensional Dirichlet distribution and Dirichlet form $\widetilde{\EE}_{\alpha,1}^{(\infty)}.$
\beg{thm}\label{weight}
Let $\gamma_{i}\geq 1,~i\geq 1.$ Denote $$\widetilde{\EE}_{\alpha,1}^{(\infty)}(f,f):=\mu_{\alpha,\alpha_{\infty}}^{(\infty)}\bigg(\sum_{i=1}^{\infty}\gamma_{i}y_i(1-|y|_{1})(\partial_{i}f)^{2}\bigg), \quad f\in\D(\widetilde{\EE}_{\alpha,1}^{(\infty)}).$$
Then the weighted super Poincar\'e inequality
$$\mu_{\alpha,\alpha_{\infty}}^{(\infty)}(f^{2})\leq r\widetilde{\EE}_{\alpha,1}^{(\infty)}(f,f)
+\beta^{(1)}(r)\mu_{\alpha,\alpha_{\infty}}^{(\infty)}(\mid f\mid)^{2},\quad f\in \mathcal{F}C^{1}(\Delta^{(\infty)})$$
holds,
 there is a positive constat $c_{n}$ such that
$$\beta^{(1)}(r)\leq  c_{n}\bigg(\frac{r}{3}\bigg)^{-[\sum_{i=1}^n 1\lor (2\aa_i)+(\aa_{\infty}-1)^+]},$$
where $n$ comes from the smallest value which satisfies
$$\frac{1}{(\sum_{i=1}^{n}\alpha_{i}+\alpha_{\infty})\inf_{i>n}\gamma_{i}}\leq r.$$
\end{thm}
\beg{proof}
$\forall r> 0,$  $n$ is the smallest value which satisfies
$$\frac{1}{(\sum_{i=1}^{n}\alpha_{i}+\alpha_{\infty})\inf_{i>n}\gamma_{i}}\leq r.$$
$\forall f\in \mathcal{F}C^{1}(\Delta^{(\infty)}),$ and the number of the variables of $f$ is denoted as $m.$
If $m\leq n,$ by the result of \cite{WZ}, we have the super Poincar\'e inequality
\beg{align*}&\mu_{\alpha,\alpha_{\infty}}^{(\infty)}(f^{2})\\
&\leq r\mu_{\alpha,\alpha_{\infty}}^{(\infty)}\bigg(\sum_{i=1}^{\infty}x_i(1-|x|_{1})(\partial_{i}f)^{2}\bigg)+\beta_{n}(r)\mu_{\alpha,\alpha_{\infty}}^{(\infty)}(\mid f\mid )^{2}\\
&\leq r\mu_{\alpha,\alpha_{\infty}}^{(\infty)}\bigg(\sum_{i=1}^{\infty}x_i\gamma_{i}(1-|x|_{1})(\partial_{i}f)^{2}\bigg)+\beta_{n}(r)\mu_{\alpha,\alpha_{\infty}}^{(\infty)}(\mid f\mid )^{2}\\
&= r\widetilde{\EE}_{\alpha,1}^{(\infty)}(f,f)
+\beta^{(1)}(r)\mu_{\alpha,\alpha_{\infty}}^{(\infty)}(\mid f\mid)^{2},\end{align*}
there is a positive constant $c_{n}$ such that
$$\beta^{(1)}(r)\leq  c_{n}\bigg(\frac{r}{3}\bigg)^{-[\sum_{i=1}^n 1\lor (2\aa_i)+(\aa_{\infty}-1)^+]}.$$
If $m> n,$ according to \cite[Theorem 1.2]{FMW}, we have the Poincar\'e inequality for $\mu_{2}^{(\infty)},$
$$\mu_{2}^{(\infty)}(f^{2})\leq \frac{1}{\sum_{i=1}^{n}\alpha_{i}+\alpha_{\infty}}\mu_{2}^{(\infty)}\bigg(\sum_{i=n+1}^{\infty}x_i(1-|x|_{1})(\partial_{i}f)^{2}\bigg)+\mu_{2}^{(\infty)}(\mid f\mid)^{2}.$$
According to \cite[Theorem 1.1]{WZ}, we have the super Poincar\'e inequality for $\mu_{1},$
$$\mu_{1}(f^{2})\leq r\mu_{1}\bigg(\sum_{i=1}^{n}x_i(1-|x|_{1})(\partial_{i}f)^{2}\bigg)+\beta_{n}(r)\mu_{1}(\mid f\mid )^{2},$$
there is a positive constant $c_{n}$ such that
 $$\beta_{n}(r)\leq c_{n}r^{-[\sum_{i=1}^n 1\lor (2\aa_i)+(\aa_{\infty}-1)^+]}.$$
Combining them together we get
\beg{align*}
&\mu_{\alpha,\alpha_{\infty}}^{(\infty)}(f^{2})=\mu_{1}(\mu_{2}^{(\infty)}(f^{2}\circ T))\\
&\leq\mu_{1}\bigg(\frac{1}{\sum_{i=1}^{n}\alpha_{i}+\alpha_{\infty}}\mu_{2}^{(\infty)}\bigg(\sum_{i=n+1}^{\infty}x_i(1-|x|_{1})(\partial_{i}(f\circ T))^{2}\bigg)\bigg)+\mu_{1}(\mu_{2}^{(\infty)}(\mid f\circ T\mid)^{2})\\
&=\mu_{1}\bigg(\frac{1}{\sum_{i=1}^{n}\alpha_{i}+\alpha_{\infty}}\mu_{2}^{(\infty)}\bigg(\sum_{i=n+1}^{\infty}x_i(1-|x|_{1})\bigg((\partial_{i}f)\circ T-\sum_{j=1}^{n}x_{j}(\partial_{j}f)\circ T\bigg)^{2} \bigg)\bigg)\\
&{\quad}+r\mu_{1}\bigg(\mu_{2}^{(\infty)}\bigg(\sum_{i=1}^{n}x_i(1-|x|_{1})(\partial_{i}(f\circ T))^{2}\bigg)\bigg)+ \beta_{n}(r)\mu_{\alpha,\alpha_{\infty}}^{(\infty)}(\mid f\mid)^{2}\\
&\leq2\mu_{1}\bigg(\frac{1}{\sum_{i=1}^{n}\alpha_{i}+\alpha_{\infty}}\mu_{2}^{(\infty)}\bigg(\sum_{i=n+1}^{\infty}x_i(1-|x|_{1})(\partial_{i}f)^{2}\circ T\bigg)\bigg)\\
&{\quad}+2\mu_{1}\bigg(\frac{1}{\sum_{i=1}^{n}\alpha_{i}+\alpha_{\infty}}\mu_{2}^{(\infty)}\bigg(\sum_{i=n+1}^{\infty}x_i(1-|x|_{1})\bigg( \sum_{j=1}^{n}x_{j}^{2}(\partial_{j}f)^{2}\circ T \bigg)\bigg)\bigg)\\
&{\quad}+r\mu_{1}\bigg(\mu_{2}^{(\infty)}\bigg(\sum_{i=1}^{n}x_i(1-|x|_{1})\bigg(\partial_{i}f\circ T\cdot\bigg(1-\sum_{i=n+1}^{\infty}x_i\bigg)\bigg)^{2}\bigg)\bigg)\\
&{\quad}+\beta_{n}(r) \mu_{\alpha,\alpha_{\infty}}^{(\infty)}(\mid f\mid)^{2}\\
&\leq2\mu_{1}\bigg(\frac{1}{\sum_{i=1}^{n}\alpha_{i}+\alpha_{\infty}}\mu_{2}^{(\infty)}\bigg(\bigg[\sum_{i=n+1}^{\infty}x_i(1-|x|_{1})(\partial_{i}f)^{2}\bigg]\circ T\bigg)\bigg)\\
&{\quad}+2\mu_{1}\bigg(\frac{1}{\sum_{i=1}^{n}\alpha_{i}+\alpha_{\infty}}\mu_{2}^{(\infty)}\bigg(\bigg[\sum_{i=1}^{n}x_i(1-|x|_{1})(\partial_{i}f)^{2}\bigg]\circ T\bigg)\bigg)\\
&{\quad}+r\mu_{1}\bigg(\mu_{2}^{(\infty)}\bigg(\bigg[\sum_{i=1}^{n}x_i(1-|x|_{1})(\partial_{i}f)^{2}\bigg]\circ T\bigg)\bigg)\\
&{\quad}+\beta_{n}(r) \mu_{\alpha,\alpha_{\infty}}^{(\infty)}(\mid f\mid)^{2}\\
&\leq\frac{2}{\sum_{i=1}^{n}\alpha_{i}+\alpha_{\infty}}\mu_{\alpha,\alpha_{\infty}}^{(\infty)}\bigg(\sum_{i=1}^{\infty}x_i(1-|x|_{1})(\partial_{i}f)^{2} \bigg)+r\mu_{\alpha,\alpha_{\infty}}^{(\infty)}\bigg(\sum_{i=1}^{n}x_i(1-|x|_{1})(\partial_{i}f)^{2}\bigg)\\
&{\quad}+\beta_{n}(r)\mu_{\alpha,\alpha_{\infty}}^{(\infty)}(\mid f\mid)^{2}\\
&\leq \frac{2}{(\sum_{i=1}^{n}\alpha_{i}+\alpha_{\infty})\inf_{i\geq 1}\gamma_{i}}\mu_{\alpha,\alpha_{\infty}}^{(\infty)}\bigg(\sum_{i=1}^{\infty}x_i\gamma_{i}(1-|x|_{1})(\partial_{i}f)^{2}\bigg)\\
&{\quad}+r\mu_{\alpha,\alpha_{\infty}}^{(\infty)}\bigg(\sum_{i=1}^{n}x_i(1-|x|_{1})(\partial_{i}f)^{2}\bigg)+\beta_{n}(r)\mu_{\alpha,\alpha_{\infty}}^{(\infty)}(\mid f\mid)^{2}\\
&\leq 3r\mu_{\alpha,\alpha_{\infty}}^{(\infty)}\bigg(\sum_{i=1}^{\infty}x_i\gamma_{i}(1-|x|_{1})(\partial_{i}f)^{2}\bigg)+\beta_{n}(r)\mu_{\alpha,\alpha_{\infty}}^{(\infty)}(\mid f\mid)^{2}.
\end{align*}
We have used
 $$\gamma_{i}\geq 1, i\geq 1$$
and let
 $$\frac{1}{\inf_{i> n}\bigg(\sum_{i=1}^{n}\alpha_{i}+\alpha_{\infty}\bigg)\gamma_{i}}\leq r$$
to get these inequalities.
So we get
$$\mu_{\alpha,\alpha_{\infty}}^{(\infty)}(f^{2})\leq r\widetilde{\EE}_{\alpha,1}^{(\infty)}(f,f)+\beta_{n}\bigg(\frac{r}{3}\bigg)\mu_{\alpha,\alpha_{\infty}}^{(\infty)}(\mid f\mid)^{2},$$
where $n$ is the smallest value which satisfies  $$\frac{1}{\inf_{i>n}(\sum_{i=1}^{n}\alpha_{i}+\alpha_{\infty})\gamma_{i}}\leq r .$$
\end{proof}

\section{The weighted super Poincar\'e  inequality for the second type Dirichlet form}
 We can obtain the weighted super Poincar\'e inequality for infinite-dimensional Dirichlet distribution and Dirichlet form $\widetilde{\EE}_{\alpha,1}^{(\infty)}$ similar to Theorem \ref{weight}.

\beg{thm}\label{weight2}
Let $\gamma_{i}\geq 1,~i\geq 1.$ Denote $$\widetilde{\EE}_{\alpha,2}^{(\infty)}(f,f)=\mu_{\alpha,\alpha_{\infty}}^{(\infty)}\bigg(\sum_{i,j=1}^{\infty}\frac{1}{1-\sum_{i>n}x_{i}}\gamma_{i}\gamma_{j}y_i(\delta_{ij}-y_{j})\partial_{i}f\partial_{j}f\bigg), \quad f\in\D(\widetilde{\EE}_{\alpha,2}^{(\infty)}).$$
Then the weighted super Poincar\'e inequality
$$\mu_{\alpha,\alpha_{\infty}}^{(\infty)}(f^{2})\leq r\widetilde{\EE}_{\alpha,2}^{(\infty)}(f,f)
+\beta^{(2)}(r)\mu_{\alpha,\alpha_{\infty}}^{(\infty)}(\mid f\mid)^{2}, \quad f\in \mathcal{F}C^{1}(\Delta^{(\infty)})$$
holds, there is a constant $c_{n}$ such that
$$\beta^{(2)}(r)\leq  c_{n}\bigg(\frac{r}{2}\bigg)^{-[\sum_{i=1}^n 1\lor (2\aa_i)+(\aa_{\infty}-1)^+]},$$
where $n$ comes from the smallest value which satisfies $$\frac{1}{(\sum_{i=1}^{n}\alpha_{i}+\alpha_{\infty})\inf_{i>1}\gamma_{i}}\leq r .$$
\end{thm}
\beg{proof}
$\forall r> 0,$  $n$ is the smallest value which satisfies
$$\frac{1}{(\sum_{i=1}^{n}\alpha_{i}+\alpha_{\infty})\inf_{i>1}\gamma_{i}}\leq r.$$
$\forall f\in \mathcal{F}C^{1}(\Delta^{(\infty)}),$ the number of the variables of $f$ denoted as $m.$
If $m\leq n,$ by the result of \cite{WZ}, we have the super Poincar\'e inequality
\beg{align*}&\mu_{\alpha,\alpha_{\infty}}^{(\infty)}(f^{2})\\
&\leq r\mu_{\alpha,\alpha_{\infty}}^{(\infty)}\bigg(\sum_{i,j=1}^{n}x_i(\delta_{ij}-x_{j})(\partial_{i}f\partial_{j}f)\bigg)+\beta_{n}(r)\mu_{\alpha,\alpha_{\infty}}^{(\infty)}(\mid f\mid )^{2}\\
&\leq r\widetilde{\EE}_{\alpha,2}^{(\infty)}(f,f)+\beta^{(2)}(r)\mu_{\alpha,\alpha_{\infty}}^{(\infty)}(\mid f\mid)^{2}
,\end{align*}
there is a positive constant $c_{n}$ such that
$$\beta^{(2)}(r)\leq  c_{n}r^{-[\sum_{i=1}^n 1\lor (2\aa_i)+(\aa_{\infty}-1)^+]}.$$
If $m> n,$ according to \cite[Proposition 3.3]{S}, we have the Poincar\'e inequality for $\mu_{2}^{(\infty)}$
$$\mu_{2}^{(\infty)}(f^{2})\leq \frac{1}{\sum_{i=n+1}^{\infty}\alpha_{i}+\alpha_{\infty}}\mu_{2}^{(\infty)}\bigg(\sum_{i,j=n+1}^{\infty}x_i(\delta_{ij}-x_{j})(\partial_{i}f\partial_{j}f)(x)\bigg)+\mu_{2}^{(\infty)}(\mid f\mid)^{2}.$$
According to \cite[Theorem 1.1]{WZ}, we have the super Poincar\'e inequality for $\mu_{1},$
$$\mu_{1}(f^{2})\leq r\mu_{1}\bigg(\sum_{i,j=1}^{n}x_i(\delta_{ij}-x_{j})(\partial_{i}f\partial_{j}f)(x)\bigg)+\beta_{n}(r)\mu_{1}(\mid f\mid )^{2},$$
there is a positive constant $c_{n}$ such that
$$\beta_{n}(r)\leq c_{n}r^{-[\sum_{i=1}^n 1\lor (2\aa_i)+(\aa_{\infty}-1)^+]}.$$
Combining them together we get
\beg{align*}
&\mu_{\alpha,\alpha_{\infty}}^{(\infty)}(f^{2})=\mu_{1}(\mu_{2}^{(\infty)}(f^{2}\circ T))\\
&\leq\mu_{1}\bigg(\frac{1}{\sum_{i=n+1}^{\infty}\alpha_{i}+\alpha_{\infty}}\mu_{2}^{(\infty)}\bigg(\sum_{i,j=n+1}^{\infty}x_i(\delta_{ij}-x_{j})\partial_{i}(f\circ T) \partial_{j}(f\circ T)\bigg)\bigg)\\
&{\quad}+\mu_{1}(\mu_{2}^{(\infty)}(\mid f\circ T\mid)^{2})\\
&=\mu_{1}\bigg(\frac{1}{\sum_{i=n+1}^{\infty}\alpha_{i}+\alpha_{\infty}}\mu_{2}^{(\infty)}\bigg(\sum_{i,j=n+1}^{\infty}x_i(\delta_{ij}-x_{j})\bigg([\partial_{i}f\circ T-\sum_{k=1}^{n}x_{k}\partial_{k}f\circ T]\\
&{\quad}\cdot[\partial_{j}f\circ T-\sum_{k=1}^{n}x_{k}\partial_{k}f\circ T]\bigg) \bigg)\bigg)
+r\mu_{1}\bigg(\mu_{2}^{(\infty)}\bigg(\sum_{i,j=1}^{n}x_i(\delta_{ij}-x_{j})\partial_{i}(f\circ T)\partial_{j}(f\circ T)\bigg)\bigg)\\
&{\quad}+\beta_{n}(r)\mu_{\alpha,\alpha_{\infty}}^{(\infty)}(\mid f\mid)^{2}\\
&\leq\mu_{1}\bigg(\frac{1}{\sum_{i=1}^{n}\alpha_{i}+\alpha_{\infty}}\mu_{2}^{(\infty)}\bigg(\bigg[\sum_{i,j=n+1}^{\infty}y_i(\delta_{ij}-y_{j})\partial_{i}f\partial_{j}f\bigg]\circ T\bigg)\bigg)\\
&{\quad}-2\mu_{1}\bigg(\frac{1}{\sum_{i=1}^{n}\alpha_{i}+\alpha_{\infty}}\mu_{2}^{(\infty)}\bigg(\sum_{i=n+1}^{\infty}\sum_{k=1}^{n}y_iy_{k}\partial_{i}f \partial_{k}f \bigg)\circ T\bigg)\bigg)\\
&{\quad}+\mu_{1}\bigg(\frac{1}{\sum_{i=1}^{n}\alpha_{i}+\alpha_{\infty}}\mu_{2}^{(\infty)}\bigg( \sum_{j=n+1}^{\infty}y_{j}\bigg(1-\sum_{j=n+1}^{\infty}y_{j}\bigg)\bigg( \sum_{k=1}^{n}x_{k}\partial_{k}f\circ T\bigg)^{2}\bigg)\bigg)\\
&{\quad}+r\mu_{1}\bigg(\mu_{2}^{(\infty)}\bigg(\sum_{i,j=1}^{n}x_i\bigg(1-\sum_{j=n+1}^{\infty}x_{j}\bigg)^{2}(\delta_{ij}-x_{j})(\partial_{i}f\partial_{j}f)\circ T\bigg)\bigg)\\
&{\quad}+\beta_{n}(r)\mu_{\alpha,\alpha_{\infty}}^{(\infty)}(\mid f\mid)^{2}\\
&\leq\frac{1}{\bigg(\sum_{i=1}^{n}\alpha_{i}+\alpha_{\infty}\bigg)}\mu_{\alpha,\alpha_{\infty}}^{(\infty)}\bigg(\sum_{i,j=1}^{\infty}\frac{1}{1-\sum_{i>n}x_{i}}x_i(\delta_{ij}-x_{j})\partial_{i}f\partial_{j}f\bigg)\\
&{\quad}+r\mu_{\alpha,\alpha_{\infty}}^{(\infty)}\bigg(\sum_{i,j=1}^{n}x_i(\delta_{ij}-x_{j})\partial_{i}f\partial_{j}f\bigg)+ \beta_{n}(r)\mu_{\alpha,\alpha_{\infty}}^{(\infty)}(\mid f\mid)^{2}\\
&\leq\frac{1}{\bigg(\sum_{i=1}^{n}\alpha_{i}+\alpha_{\infty}\bigg)\inf_{i> n}\gamma_{i}}\mu_{\alpha,\alpha_{\infty}}^{(\infty)}\bigg(\sum_{i,j=1}^{\infty}\frac{1}{1-\sum_{i>n}x_{i}}\gamma_{i}\gamma_{j}x_i(\delta_{ij}-x_{j})\partial_{i}f\partial_{j}f\bigg)\\
&{\quad}+r\mu_{\alpha,\alpha_{\infty}}^{(\infty)}\bigg(\sum_{i,j=1}^{n}x_i(\delta_{ij}-x_{j})\partial_{i}f\partial_{j}f\bigg)+ \beta_{n}(r)\mu_{\alpha,\alpha_{\infty}}^{(\infty)}(\mid f\mid)^{2}\\
&\leq 2r\mu_{\alpha,\alpha_{\infty}}^{(\infty)}\bigg(\sum_{i,j=1}^{\infty}\frac{1}{1-\sum_{i>n}x_{i}}\gamma_{i}\gamma_{j}y_i(\delta_{ij}-y_{j})\partial_{i}f\partial_{j}f\bigg)
+\beta_{n}(r)\mu_{\alpha,\alpha_{\infty}}^{(\infty)}(\mid f\mid)^{2}.
\end{align*}
We have used $\gamma_{i}\geq 1,~ i\geq 1$ and $$\frac{1}{\bigg(\sum_{i=1}^{n}\alpha_{i}+\alpha_{\infty}\bigg)\inf_{i> n}\gamma_{i}}\leq r$$
to get these inequalities.
So we get
$$\mu_{\alpha,\alpha_{\infty}}^{(\infty)}(f^{2})\leq r\widetilde{\EE}_{\alpha,2}^{(\infty)}(f,f)
+\beta_{n}(r)\mu_{\alpha,\alpha_{\infty}}^{(\infty)}(\mid f\mid)^{2},$$
there is a positive constant $c_{n}$ such that
$$\beta_{n}(r)\leq c_{n}\bigg(\frac{r}{2}\bigg)^{-[\sum_{i=1}^n 1\lor (2\aa_i)+(\aa_{\infty}-1)^+]},$$
where $n$ is the smallest value which satisfies  $$\frac{1}{\inf_{i>n}(\sum_{i=1}^{n}\alpha_{i}+\alpha_{\infty})\gamma_{i}}\leq r .$$
\end{proof}

\section*{Acknowledgement} The author would like to thank Professor Feng-Yu Wang
for helpful guidance and other's suggestions.

\section*{Funding}
The work is supported in part by NNSFC (11771326).

\section*{Competing interests}
The authors declare that they have no competing interests.

\beg{thebibliography}{99}
\bibitem{ACM} A. Albanese, M. Campiti, E. Mangino, \emph{Regularity properties of semigroups generated by some Fleming-type operators,} J. Math. Anal. Appl. 335 (2007) 1259--1273.

\bibitem{BR} J. Bakosi, J. R. Ristorcelli, \emph{A stochastic diffusion process for the Dirichlet distribution,} International J. Stoch. Anal. 2013, Article ID 842981, 7 pages.

\bibitem{FMW} S. Feng, L. Miclo, F.-Y. Wang, \emph{ Poincar\'e inequality for Dirichlet distributions and infinite-dimensional generalizations,} Lat. Am. J. Probab. Math. Stat. 14(2017), 361--380.

 \bibitem{FW07}  S. Feng, F.-Y. Wang,  \emph{A class of infinite-dimensional diffusion processes with connection to population genetics,} J. Appl. Probab. 44(2007), 938--949.

 \bibitem{Gross} L. Gross, \emph{ Logarithmic Sobolev inequalities
and contractivity properties of semigroups,}  in
$``$Dirichlet Forms'',  Lecture Notes in
Math.  1563 (Springer, Berlin), pp. 54--88.

\bibitem{Jac01} M. Jacobsen, \emph{Examples of multivariate diffusions: time-reversibility; a Cox-Ingersoll-Ross type process,} Department of Theoretical Statistics, Preprint 6, University of Copenhagen, 2001.

   \bibitem{P1} N. L. Johnson, \emph{An approximation to the multinomial distribution, some properties and applications,} Biometrika, 47(1960), 93--102.

   \bibitem{Mi1} L. Miclo, \emph{About projections of logarithmic Sobolev inequalities,}  Lecture Notes in Math. 1801 (J. Az\'ema, M. \'Emery, M. Ledoux, M. Yor Eds), pp. 201--221, 2003, Springer.

 \bibitem{Mi2} L. Miclo, \emph{Sur l'in\'egalit\'e de Sobolev logarithmique des op\'erateurs de Laguerre \`{a} petit param\'etre, } Lecture Notes in Math. 1801 (J. Az\'ema, M. \'Emery, M. Ledoux, M. Yor Eds), pp. 222--229, 2003, Springer.

  \bibitem{S} W. Stannat, \emph{On validity of the log-Sobolev
      inequality for symmetric Fleming-Viot operators,}  Ann.
      Probab.
      28(2000), 667--684.

\bibitem{W00a}
F.-Y. Wang, \emph{ Functional inequalities for empty essential
spectrum, }  J. Funct. Anal. 170(2000), 219--245.

\bibitem{W00b} F.-Y. Wang, \emph{ Functional inequalities, semigroup properties
and spectrum estimates,}   Infin. Dimens. Anal. Quant. Probab.
Relat. Topics 3(2000), 263--295.

 \bibitem{Wbook} F.-Y. Wang, \emph{Functional Inequalities, Markov Semigroups and Spectral Theory,}   Science Press 2005.

\bibitem{WZ} F.-Y. Wang, W.W. Zhang, \emph{Nash inequality for diffusion processes associated with Dirichlet distributions,} 2018.

\end{thebibliography}
\end{document}